\documentclass{ifacconf}

\usepackage{graphicx}      % include this line if your document contains figures
\usepackage{natbib}        % required for bibliography
\usepackage{amsmath,amssymb,mathrsfs}

\usepackage{tikz}
\usetikzlibrary{decorations.pathreplacing}
\usetikzlibrary{patterns}

\newtheorem{theorem}{Theorem}

\newtheorem{remark}{Remark}

\DeclareMathOperator*{\argmax}{arg\,max}

\def\R{\mathbb{R}}
\def\N{\mathbb{N}}

\def\t{\tau}
\def\S{\mathrm{S}}
\def\CC{\mathscr{C}}

\def\OCPI{{\bf (OCP)}}
\def\OSCPI{{\bf (OSDCP)}}

\begin{document}
\begin{frontmatter}

\title{Pontryagin maximum principle for optimal sampled-data control problems.}

%\title{Pontryagin maximum principle for optimal sampled-data control problems.\thanksref{footnoteinfo}} 

%\thanks[footnoteinfo]{Sponsor and financial support acknowledgment
%goes here. Paper titles should be written in uppercase and lowercase
%letters, not all uppercase.}

\author[First]{Lo\"ic Bourdin} 
\author[Second]{Emmanuel Tr\'elat}

\address[First]{Universit\'e de Limoges, Institut de recherche XLIM, CNRS UMR 7252, D\'epartement de Math\'ematiques et d'Informatique, Limoges, France (e-mail: loic.bourdin@unilim.fr).}
\address[Second]{Sorbonne Universit\'es, UPMC Univ Paris 06, CNRS UMR 7598, Laboratoire Jacques-Louis Lions, Institut Universitaire de France, F-75005, Paris, France (e-mail: emmanuel.trelat@upmc.fr)}

\begin{abstract}                % Abstract of not more than 250 words.
In this short communication, we first recall a version of the Pontryagin maximum principle for general finite-dimensional nonlinear optimal sampled-data control problems. This result was recently obtained in [L.~Bourdin and E.~Tr\'elat, \textit{Optimal sampled-data control, and generalizations on time scales}, arXiv:1501.07361, 2015]. Then we discuss the maximization condition for optimal sampled-data controls that can be seen as an average of the weak maximization condition stated in the classical Pontryagin maximum principle for optimal (permanent) controls. Finally, applying this theorem, we solve a linear-quadratic example based on the classical parking problem.
\end{abstract}

\end{frontmatter}
%===============================================================================

\section{Introduction}
Optimal control theory is concerned with the analysis of controlled dynamical systems, where one aims at steering such a system from a given configuration to some desired target by minimizing some criterion. The Pontryagin maximum principle (in short, PMP), established at the end of the 50's for general finite-dimensional nonlinear continuous-time dynamics (see \cite{pont}, and see \cite{gamk} for the history of this discovery), is certainly the milestone of the classical optimal control theory. It provides a first-order necessary condition for optimality, by asserting that any optimal trajectory must be the projection of an extremal. The PMP then reduces the search of optimal trajectories to a boundary value problem posed on extremals. Optimal control theory, and in particular the PMP, has an immense field of applications in various domains, and it is not our aim here to list them.

%\medskip

We speak of a {\em purely continuous-time optimal control problem}, when both the state $q$ and the control $u$ evolve continuously in time, and the control system under consideration has the form
$$
\dot q(t) = f(t,q(t),u(t)), \; \text{for a.e. } t \in \R^+,
$$
where $q(t) \in \R^n$ and $u(t) \in \Omega \subset \R^m$. Such models assume that the control is permanent, that is, the value of $u(t)$ can be chosen at each time $t \in \R^+$. We refer the reader to textbooks on continuous optimal control theory such as \cite{agrach,Bon-Chy03,trel2,bres,brys,BulloLewis,hest,Jurdjevic,lee,pont,Schattler,seth,trel} for many examples of theoretical or practical applications.

% \medskip

We speak of a {\em purely discrete-time optimal control problem}, when both the state $q$ and the control $u$ evolve in a discrete way in time, and the control system under consideration has the form
$$
q_{k+1}-q_k = f(k,q_k,u_k), \; k \in \N,
$$
where $q_k \in \R^n$ and $u_k \in \Omega \subset \R^m$. As in the continuous case, such models assume that the control is permanent, that is, the value of $u_k$ can be chosen at each time $k \in \N$. A version of the PMP for such discrete-time control systems has been established in \cite{Halkin,holt2,holt} under appropriate convexity assumptions. The considerable development of the discrete-time control theory was in particular motivated by the need of considering digital systems or discrete approximations in numerical simulations of differential control systems (see the textbooks \cite{bolt,cano,mord,seth}).
It can be noted that some early works devoted to the discrete-time PMP (like \cite{fan}) are mathematically incorrect. Some counterexamples were provided in \cite{bolt} (see also \cite{mord}), showing that, as is now well known, the exact analogous of the continuous-time PMP does not hold at the discrete level. More precisely, the maximization condition of the continuous-time PMP cannot be expected to hold in general in the discrete-time case. Nevertheless, a weaker condition can be derived, in terms of nonpositive gradient condition (see Theorem 42.1 in \cite{bolt}).

% \medskip

We speak of an {\em optimal sampled-data control problem}, when the state $q$ evolves continuously in time, whereas the control $u$ evolves in a discrete way in time. This hybrid situation is often considered in practice for problems in which the evolution of the state is very quick (and thus can be considered continuous) with respect to that of the control. We often speak, in that case, of {\em digital control}. This refers to a situation where, due for instance to hardware limitations or to technical difficulties, the value $u(t)$ of the control can be chosen only at times $t=kT$, where $T>0$ is fixed and $k \in \N$. This means that, once the value $u(kT)$ is fixed, $u(t)$ remains constant over the time interval $[kT,(k+1)T)$. Hence the trajectory $q$ evolves according to
%{\small
$$
\dot q(t) = f(t,q(t),u(kT)), \; \text{for a.e. } t \in [kT,(k+1)T),\; k\in\N.
$$
%}
In other words, this {\em sample-and-hold} procedure consists of ``freezing" the value of $u$ at each \textit{controlling time} $t=kT$ on the corresponding \textit{sampling time interval} $[kT,(k+1)T)$, where $T$ is called the \textit{sampling period}. In this situation, the control of the system is clearly nonpermanent.

% \medskip

To the best of our knowledge, the classical optimal control theory does not treat general nonlinear optimal sampled-data control problems, but concerns either purely continuous-time, or purely discrete-time optimal (permanent) control problems. In \cite{bourdin-trelat-pontryagin2} we provided a version of the PMP that can be applied to general nonlinear optimal sampled-data control problems.\footnote{Actually we established in \cite{bourdin-trelat-pontryagin2} a PMP in the much more general framework of {\em time scales}, which unifies and extends continuous-time and discrete-time issues. But it is not our aim here to enunciate this result in its whole generality.}

% \medskip

In this short communication, we first recall in Section~\ref{section1} the above mentioned PMP. Then a discussion is provided concerning the maximization condition for optimal sampled-data controls that can be seen as an average of the weak maximization condition stated in the classical PMP for optimal (permanent) controls. Finally, in Section~\ref{section2}, we solve a linear-quadratic example based on the classical parking problem.

\section{Main result}\label{section1}

Let $m$, $n$ and $j$ be nonzero integers. In the sequel, we denote by $\langle \cdot , \cdot \rangle_n$ the classical scalar product in $\R^n$. Let $T>0$ be an arbitrary sampling period. In what follows, for any real number $t$, we denote by $E(t)$ the integer part of $t$, defined as the unique integer such that $E(t)\leq t< E(t)+1$. Note that $k = E(t/T)$ whenever $kT\leq t<(k+1)T$. 

In this section, we are interested in the general nonlinear optimal sampled-data control problem given by
{\small
\begin{equation*}
\OSCPI \; \left\{\begin{split}
& \min \int_0^{t_f} f^0 ( \t, q(\t), u( kT ) ) \, d\t , \quad \textrm{with}\ k=E(\t/T)  , \\
& \dot q(t) = f ( t,q(t), u( kT ) ), \quad \textrm{with}\ k=E(t/T)  , \\[5pt]
& u(kT) \in \Omega , \\[5pt]
& g(q(0),q(t_f)) \in \S .
\end{split}\right.
\end{equation*}}Here, $f: \R\times \R^{n} \times \R^{m}  \rightarrow \R^{n}$, $f^0: \R\times \R^{n} \times \R^{m}  \rightarrow \R$ and $g: \R^n \times \R^n \rightarrow \R^j$ are of class $\CC^1$,
and $\Omega$ (resp., $\S$) is a non-empty closed convex subset of $\R^m$ (resp., of $\R^j$).
The final time $t_f \geq 0$ can be fixed or not.

%Note that, under appropriate (usual) compactness and convexity assumptions, the optimal control problem $\OSCPI$ has at least one solution (see Theorem~\ref{propexistence} in Section~\ref{subsectionOSCP}).

Recall that $g$ is said to be submersive at a point $(q_1,q_2) \in \R^n \times \R^n$ if the differential of $g$ at this point is surjective. We define as usual the Hamiltonian $H:\R\times\R^n \times \R^n \times \R \times \R^m\rightarrow \R$ by
$$ H(t,q,p,p^0,u) = \langle p, f(t,q,u) \rangle_{n} + p^0 f^0(t,q,u) . $$ 

\subsection{Statement}

In \cite{bourdin-trelat-pontryagin2} we proved the following theorem.

\begin{theorem}[PMP for $\OSCPI$]\label{thmmainintro}
If a trajectory $q$, defined on $[0,t_f]$ and associated with a sampled-data control $u$, is an optimal solution of $\OSCPI$, then there exists a nontrivial couple $(p,p^0)$, where $p : [0,t_f] \rightarrow \R^n$ is an absolutely continuous mapping (called adjoint vector) and $p^0 \leq 0$, such that the following conditions hold:
\begin{itemize}
\item \textbf{Extremal equations}:
$$ \dot q(t) = \partial_p H (t,q(t),p(t),p^0,u(kT)), $$
$$ \dot p(t) = -\partial_q H (t,q(t),p(t),p^0,u(kT)) ,$$
for almost every $t\in[0,t_f)$, with $k=E(t/T)$. \\
\item \textbf{Maximization condition}:\\
For every controlling time $kT \in [0,t_f)$ such that $(k+1)T \leq  t_f$, we have
\begin{multline}\label{secondconditionintro}
\Big\langle \dfrac{1}{T} \int_{kT}^{(k+1)T} \partial_u H (\t,q(\t),p(\t),p^0,u(kT)) \; d\t  \\  , \; y-u(kT) \Big\rangle_{m}  \leq 0,
\end{multline}
for every $y \in \Omega$. In the case where $kT \in [0,t_f)$ with $(k+1)T >  t_f$, the above maximization condition is still valid provided $\frac{1}{T}$ is replaced with $\frac{1}{t_f - kT}$ and $(k+1)T$ is replaced with $t_f$. \\
\item \textbf{Transversality conditions on the adjoint vector}:\\
If $g$ is submersive at $(q(0),q(t_f))$, then the nontrivial couple $(p,p^0)$ can be selected to satisfy
$$ p(0) = - \partial_1 g (q(0),q(t_f))^\top  \psi, $$
$$ p(t_f) =  \partial_2 g (q(0),q(t_f)) ^\top  \psi, $$
where $-\psi$ belongs to the orthogonal of $\S$ at the point $g (q(0),q(t_f)) \in \S$. \\
\item \textbf{Transversality condition on the final time}:\\
If the final time is left free in the optimal sampled-data control problem $\OSCPI$ and if $t_f>0$, then the nontrivial couple $(p,p^0)$ can be moreover selected to satisfy
\begin{equation*}
H(t_f, q(t_f), p(t_f),p^0,u(k_f T) ) = 0,
\end{equation*}
where $k_f =E(t_f/T)$ whenever $t_f \notin \N T$, and $k_f=E(t_f/T)-1$ whenever $t_f \in \N T$.
\end{itemize}
\end{theorem}

The maximization condition~\eqref{secondconditionintro}, which is satisfied for every $y \in \Omega$, gives a necessary condition allowing to compute $u(kT)$ in general, and this, for all controlling times $kT \in [0,t_f)$. We will solve in Section~\ref{section2} an example of optimal sampled-data control problem, and show how these computations can be done in a simple way.

\begin{remark}\label{remarknormaliser}
As is well known, the nontrivial couple $(p,p^0)$ of Theorem~\ref{thmmainintro}, which is a Lagrange multiplier, is defined up to a multiplicative scalar. Defining as usual an \textit{extremal} as a quadruple $(q,p,p^0,u)$ solution of the extremal equations, an extremal is said to be \textit{normal} whenever $p^0\neq 0$ and \textit{abnormal} whenever $p^0=0$. In the normal case $p^0\neq 0$, it is usual to normalize the Lagrange multiplier so that $p^0=-1$.
\end{remark}

\begin{remark}\label{remarkconditionsterminales}
Let us describe some typical situations of terminal conditions $g(q(0),q(t_f)) \in \S$ in $\OSCPI$, and of the corresponding transversality conditions on the adjoint vector.
\begin{itemize}
\item If the initial and final points are fixed in $\OSCPI$, that is, if we impose $q(0) = q_0$ and $q(t_f) = q_f$, then $j=2n$, $g(q_1,q_2) = (q_1,q_2)$ and $\S = \{ q_0 \} \times \{ q_f \}$. In that case, the transversality conditions on the adjoint vector give no additional information.
\item If the initial point is fixed, that is, if we impose $q(0) = q_0$, and if the final point is left free in $\OSCPI$, then $j=n$, $g(q_1,q_2) = q_1$ and $\S = \{ q_0 \} $. In that case, the transversality conditions on the adjoint vector imply that $p(t_f) = 0$. Moreover, we have $p^0 \neq 0$\footnote{Indeed, if $p^0 =0$, then the adjoint vector $p$ is trivial from the extremal equation and from the final condition $p(t_f)=0$. This leads to a contradiction since the couple $(p,p^0)$ has to be nontrivial.} and we can normalize the Lagrange multiplier so that $p^0=-1$ (see Remark~\ref{remarknormaliser}).
\item If the periodic condition $q(0)=q(t_f)$ is imposed in $\OSCPI$, then $j=n$, $g(q_1,q_2) = q_1- q_2$ and $\S = \{ 0 \}$. In that case, the transversality conditions on the adjoint vector yield that $p(0) = p(t_f)$.
\end{itemize}
We stress that, in all examples above, the function $g$ is indeed a submersion.
\end{remark}

\begin{remark}
In \cite{bourdin-trelat-pontryagin2} we also provided a result stating the existence of optimal solutions for $\OSCPI$, under some appropriate compactness and convexity assumptions. Actually, if the existence of solutions is stated, the necessary conditions provided in Theorem~\ref{thmmainintro} may prove the uniqueness of the optimal solution. 
\end{remark}

\subsection{Averaging of the classical weak maximization condition}\label{section12}

Let us compare the maximization condition~\eqref{secondconditionintro} with respect to that of the classical PMP. Let us consider the following general nonlinear optimal (permanent) control problem
\begin{equation*}
\OCPI \; \left\{\begin{split}
& \min \int_0^{t_f} f^0 ( \t, q(\t), u( \t ) ) \, d\t ,  \\
& \dot q(t) = f ( t,q(t), u( t ) ),\\[5pt]
& u(t) \in \Omega , \\[5pt]
& g(q(0),q(t_f)) \in \S .
\end{split}\right.
\end{equation*}
In the sequel, we will denote by $u^*$ an optimal (permanent) control. In the case of $\OCPI$, the statement of the classical PMP coincides with that of Theorem~\ref{thmmainintro}, except the maximization condition~\eqref{secondconditionintro}.\footnote{Actually the transversality condition on the final time is slightly different. Precisely, if the final time is left free in the optimal control problem $\OCPI$ and if $t_f>0$, then the nontrivial couple $(p,p^0)$ can be selected such that the function $t \mapsto H(t,q(t),p(t),p^0,u^*(t))$ is equal almost everywhere to a continuous function vanishing at $t=t_f$.} Indeed, the maximization condition in the classical PMP is celebrated to be given by
\begin{equation}\label{eqstrongmaxcondition}
u^*(t) \in \argmax_{y \in \Omega} H((t,q(t),p(t),p^0,y)  , 
\end{equation}
for a.e. $t \in [0,t_f)$. Note that \eqref{eqstrongmaxcondition} can be directly weakened as follows:
\begin{equation}\label{eqweakmaxcondition}
 \left\langle \partial_u H (t,q(t),p(t),p^0,u^*(t)) , y-u^*(t) \right\rangle_{m}  \leq 0, 
\end{equation}
for every $y \in \Omega$ and for a.e. $t \in [0,t_f)$. If the classical PMP is stated with the nonpositive gradient condition~\eqref{eqweakmaxcondition}, the literature speaks of \textit{weak formulation of the classical PMP}.\footnote{As mentioned in the introduction, only the weak formulation of the classical PMP can be extended to the discrete case. To extend the strong formulation of the classical PMP to the discrete case, one has to consider additional convexity assumptions on the dynamics, see Remark~\ref{hamiltonianconcave} or \cite{Halkin,holt2,holt} for example.}

\medskip

It is worth to emphasize that the maximization condition~\eqref{secondconditionintro} given in Theorem~\ref{thmmainintro} can be seen as an average of the weak maximization condition~\eqref{eqweakmaxcondition} given in the classical PMP. For this reason we speak of \textit{nonpositive average gradient condition}.

\begin{remark}\label{hamiltonianconcave}
In the case where the Hamiltonian $H$ is concave in $u$, the strong and the weak formulations of the classical PMP are obviously equivalent. In a similar way, if $H$ is concave in $u$, note that the maximization condition~\eqref{secondconditionintro} in Theorem~\ref{thmmainintro} can be written as 
$$ u(kT) \in \argmax_{y \in \Omega} \,  \dfrac{1}{T} \int_{kT}^{(k+1)T} H (\t,q(\t),p(\t),p^0,y) \; d\t ,$$
for all controlling times $kT \in [0,t_f)$. In the case where $(k+1)T >  t_f$, the above maximization condition is still valid provided $\frac{1}{T}$ is replaced with $\frac{1}{t_f - kT}$ and $(k+1)T$ is replaced with $t_f$. In that case we speak of \textit{pointwise maximization of the average Hamiltonian}.
\end{remark}

\section{The parking problem}\label{section2}

In this section, we consider the classical double integrator
$$ \ddot{q} = u , \quad u \in [-1,1], $$
which can represent a car with position $q \in \R$ and with bounded acceleration $u$ acting as the control. Let us study the classical problem of parking the car at the origin, from an initial position $M > 0$ and with a fixed final time $t_f > 0$, minimizing the energy 
$$ \int_0^{t_f} u^2 \, d\t. $$
In the sequel we first give some recalls on the classical permanent control case (solved with the help of the classical PMP). Then we solve the sampled-data control case with the help of Theorem~\ref{thmmainintro} and compare the two situations.

\subsection{Recalls on the permanent control case}\label{section21}
The above optimal control problem, in the permanent control case, can be summarized as follows:
\begin{equation*}
\left\{\begin{split}
& \min \int_0^{t_f} u(\t)^2  \, d\t ,  \\[5pt]
&  \left( \begin{array}{c}
\dot{q}_1(t) \\ \dot{q}_2(t)
\end{array} \right)  = \left( \begin{array}{c}
q_2(t) \\ u(t)
\end{array} \right), \\[5pt]
& u(t) \in [-1,1] , \\[5pt]
& \left( \begin{array}{c}
q_1(0) \\ q_2(0)
\end{array} \right) =  \left( \begin{array}{c}
M \\ 0
\end{array} \right), \quad \left( \begin{array}{c}
q_1(t_f) \\ q_2(t_f)
\end{array} \right) =  \left( \begin{array}{c}
0 \\ 0
\end{array} \right).
\end{split}\right.
\end{equation*}
In the sequel we assume that $t_f^2 > 4M$ in order to ensure the existence of a solution. 

\begin{figure}[h]
\begin{center}
\begin{tikzpicture}[scale=0.8]
		\draw[->] (-0.5,0)--(10,0);
		\draw[->] (0,-2)--(0,2);
		\draw [red,domain=0:3*((3-sqrt(3))/2)] plot (\x,{-1});
		\draw [red,domain=(9-3*((3-sqrt(3))/2)):9] plot (\x,{1});		
		\draw [red,domain=3*((3-sqrt(3))/2):(9-3*((3-sqrt(3))/2))] plot (\x,{  (2*\x-9)/(3*sqrt(3))      });		
		\node at (-0.3,-0.35) {$0$};
	    \node at (9,0) {$|$};
		\node at (9,-0.5) {$t_f$};
	    \node at (1.9,0) {$|$};
		\node at (1.9,-0.5) {$t_1$};	
	    \node at (7.09,0) {$|$};
		\node at (7.09,-0.5) {$t_f-t_1$};	
		\node at (6,1) {\textcolor{red}{$u^*$}};	
	    \node at (0,1) {$-$};
	    \node at (0,-1) {$-$};
\end{tikzpicture}
\caption{Optimal (permanent) control, if $4M < t_f^2 < 6M$}\label{Fig1}
\end{center}
\end{figure}
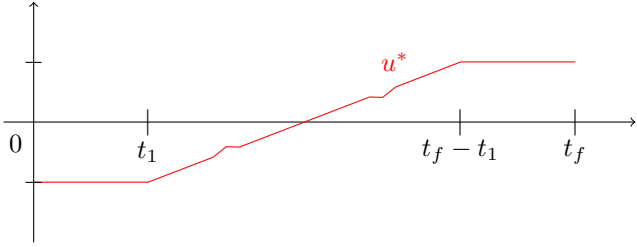

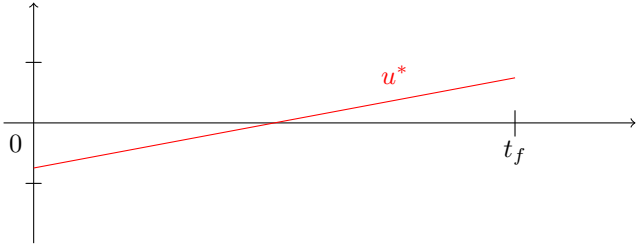
\begin{figure}[h]
\begin{center}
\begin{tikzpicture}[scale=0.8]
		\draw[->] (-0.5,0)--(10,0);
		\draw[->] (0,-2)--(0,2);
		\draw [red,domain=0:8] plot (\x,{ (12*(2*\x-8))/(16*8)});	
		\node at (-0.3,-0.35) {$0$};
	    \node at (8,0) {$|$};
		\node at (8,-0.5) {$t_f$};
		\node at (6,0.8) {\textcolor{red}{$u^*$}};	
	    \node at (0,1) {$-$};
	    \node at (0,-1) {$-$};
\end{tikzpicture}
\caption{Optimal (permanent) control, if $6M \leq  t_f^2 $}\label{Fig2}
\end{center}
\end{figure}

From the classical PMP, one can prove that, if $4M < t_f^2 < 6M$, the optimal (permanent) control $u^*$ is given by
$$ u^* (t) = \left\lbrace 
\begin{array}{lcl}
-1 & \text{if} & 0 \leq t \leq t_1, \\ \\
\dfrac{2t-t_f}{\sqrt{3(t_f^2 - 4M)}} & \text{if} & t_1 \leq t \leq t_f - t_1 , \\ \\
1 & \text{if} & t_1 \leq t \leq t_f,
\end{array}
\right. $$
where $ t_1 = \frac{1}{2} (t_f - \sqrt{3(t_f^2 - 4M)}) < \frac{t_f}{2}$, see Figure~\ref{Fig1}.
If $6M \leq t_f^2$, one can prove that the optimal (permanent) control $u^*$ is given by
$$ u^*(t) = \dfrac{6M}{t_f^3	} (2t - t_f), \quad t \in [0,t_f],$$
see Figure~\ref{Fig2}.

\subsection{The sampled-data control case}\label{sectionexamplesampled}
In this section, we consider the corresponding optimal sampled-data control problem given by
\begin{equation*}
\left\{\begin{split}
& \min \int_0^{t_f} u(kT)^2  \, d\t , \quad \textrm{with}\ k=E(\t/T)  ,  \\[5pt]
&  \left( \begin{array}{c}
\dot{q}_1(t) \\ \dot{q}_2(t)
\end{array} \right)  = \left( \begin{array}{c}
q_2(t) \\ u(kT)
\end{array} \right), \quad \textrm{with}\ k=E(t/T)  ,  \\[5pt]
& u(kT) \in [-1,1] , \\[5pt]
& \left( \begin{array}{c}
q_1(0) \\ q_2(0)
\end{array} \right) =  \left( \begin{array}{c}
M \\ 0
\end{array} \right), \quad \left( \begin{array}{c}
q_1(t_f) \\ q_2(t_f)
\end{array} \right) =  \left( \begin{array}{c}
0 \\ 0
\end{array} \right),
\end{split}\right.
\end{equation*}
where $T > 0$ is a fixed sampling period. In order to avoid the case where a controlling time $kT$ is such that $(k+1)T > t_f$ and in order to simplify the redaction, we assume that $t_f = KT$ for some $K \in \N^*$.

Let us apply Theorem~\ref{thmmainintro} in the normal case $p^0 = -1$. From the extremal equations, the adjoint vector $p=(p_1 \; p_2)^\top$ is such that $p_1$ is constant and $p_2 (t) = p_1(t_f - t)+p_2 (t_f)$ is affine. The maximization condition~\eqref{secondconditionintro} provides
$$ \frac{1}{T} (y-u(kT)) \int_{kT}^{(k+1)T} p_2(\t) - 2 u(kT) \; d\t \leq 0, $$
that is
$$ (y-u(kT)) \left[  - 2 u(kT) +p_1 \left( t_f - kT - \dfrac{T}{2} \right) + p_2 (t_f) \right] \leq 0,  $$
for all $k=0,\ldots,K-1$ and all $y \in [-1,1]$. Let us write this maximization condition as
$$ (y-u(kT)) \Gamma_k ( u(kT) ) \leq 0,  $$
for all $k=0,\ldots,K-1$ and all $y \in [-1,1]$, where $\Gamma_k : [-1,1] \to \R$ is a decreasing affine function. It clearly follows that
\begin{itemize}
\item if $\Gamma_k (-1) < 0$, then $u(kT) = -1$;
\item if $\Gamma_k (1) > 0$, then $u(kT) = 1$;
\item if $\Gamma_k (-1) > 0$ and $\Gamma_k (1) < 0$, then $u(kT)$ is the unique solution of $\Gamma_k ( x ) = 0$ given by
$$ u(kT) = \dfrac{1}{2} \left[ p_1 \left(t_f - kT - \dfrac{T}{2} \right) + p_2 (t_f)  \right] . $$
\end{itemize}
Hence, for each couple $(p_1,p_2(t_f))$, the above method allows to compute explicitly the associated values $u(kT)$ for all $k=0,\ldots,K-1$. Unfortunately, the transversality conditions on the adjoint vector do not provide any additional information on the values of $p_1$ and $p_2 (t_f)$, see Remark~\ref{remarkconditionsterminales}. As a consequence, and as usual, we proceed to a numerical shooting method on the application 
$$ (p_1,p_2(t_f)) \longmapsto (q_1(t_f),q_2(t_f)) $$ 
in order to guarantee the final constraints $q_1(t_f) = q_2 (t_f) = 0$.\footnote{In order to initiate the shooting method, we take the values of $p_1$ and $p_2 (t_f)$ from the classical permanent control case, see Section~\ref{section21}.}

Finally we obtain the following numerical results. The values $u(kT)$ are represented with blue crosses and the red curve corresponds to the optimal (permanent) control $u^*$ obtained in Section~\ref{section21}.
\begin{itemize}
\item With $M=2$, $t_f=3$ (in the case $4M < t_f^2 < 6M$) and for $T=1$, $T=0.5$, $T=0.1$ and $T=0.01$, we obtain:
\begin{center}
\includegraphics[scale=0.18]{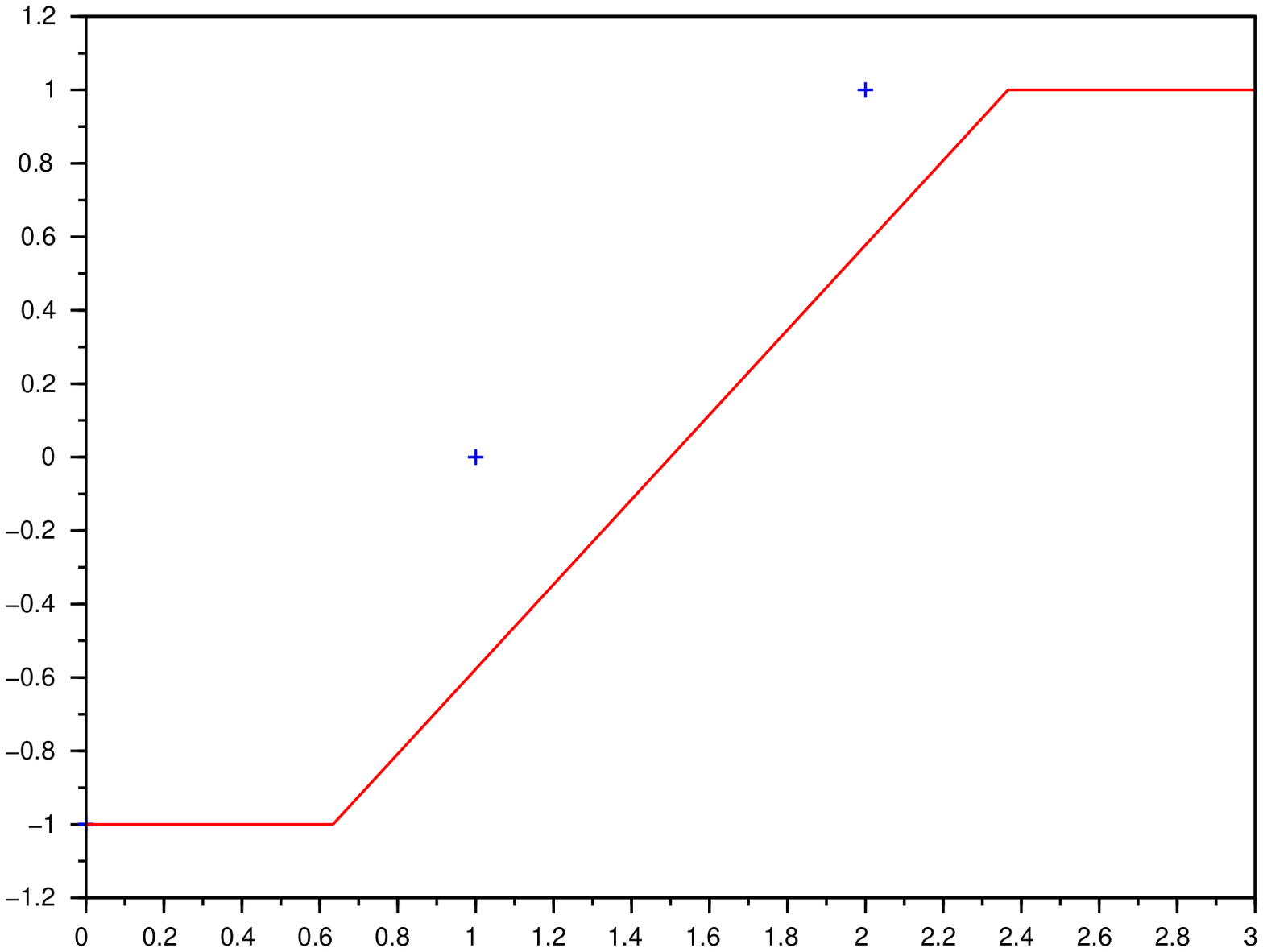}  \includegraphics[scale=0.18]{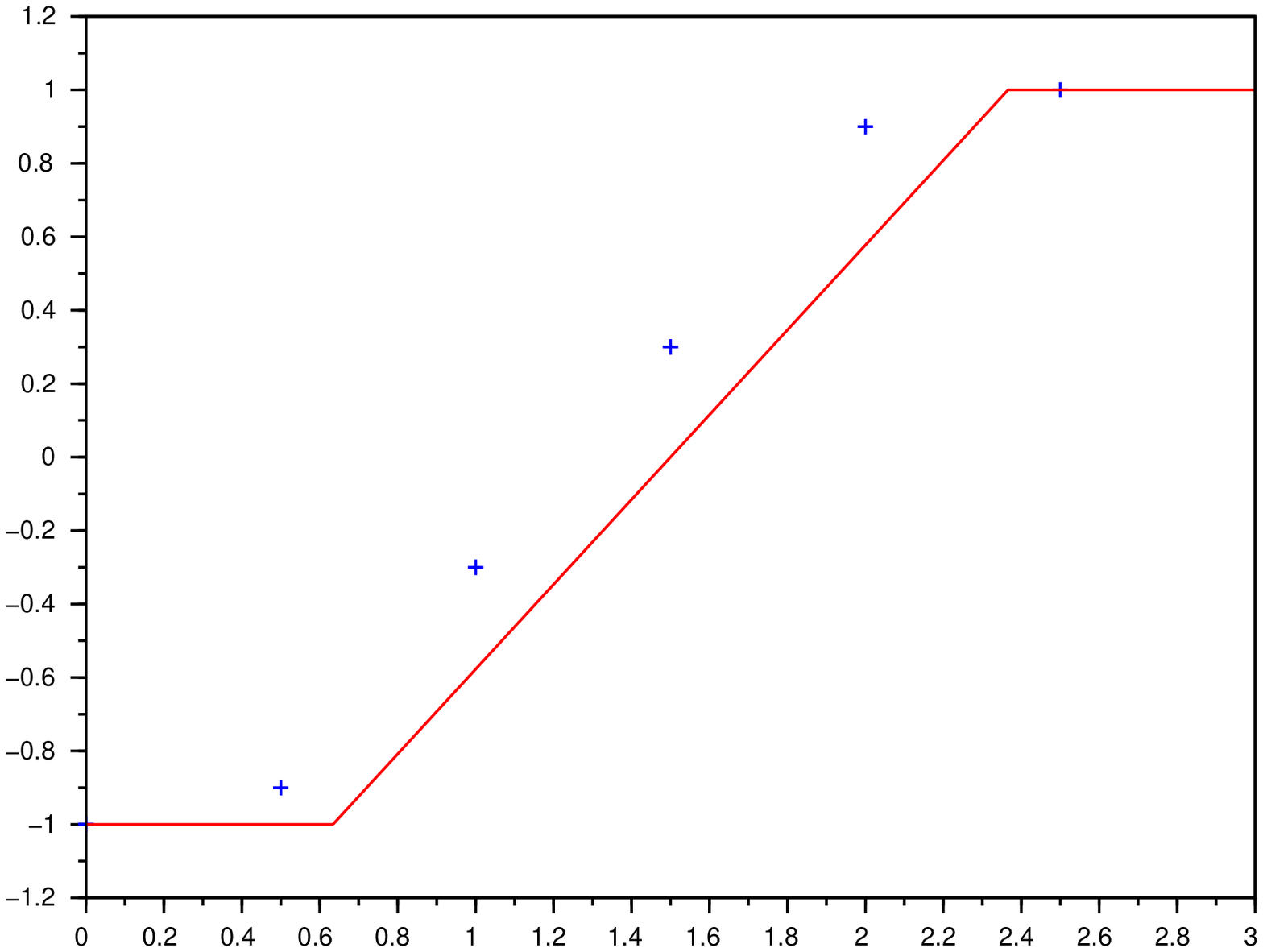} 
\includegraphics[scale=0.18]{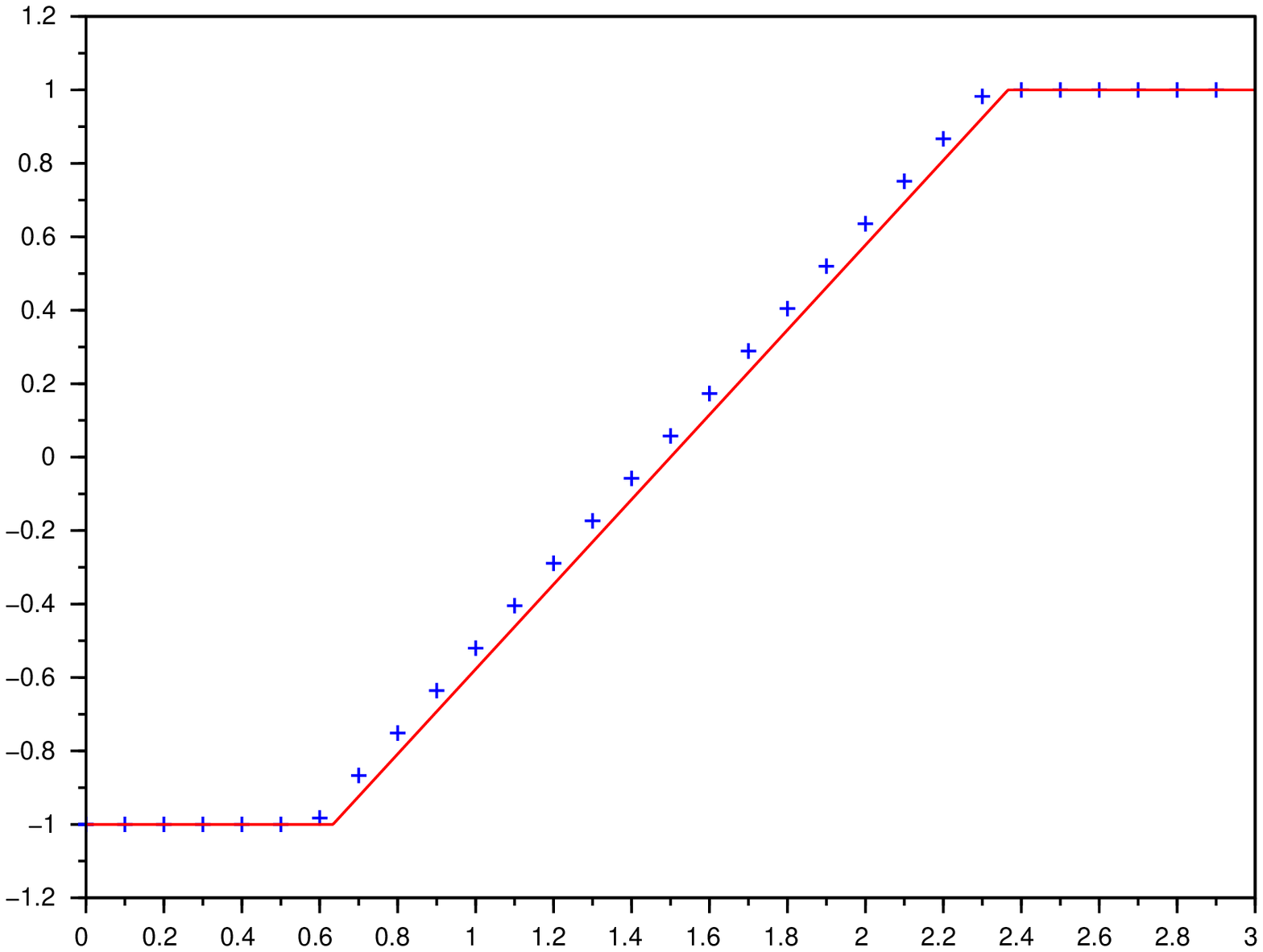}  \includegraphics[scale=0.18]{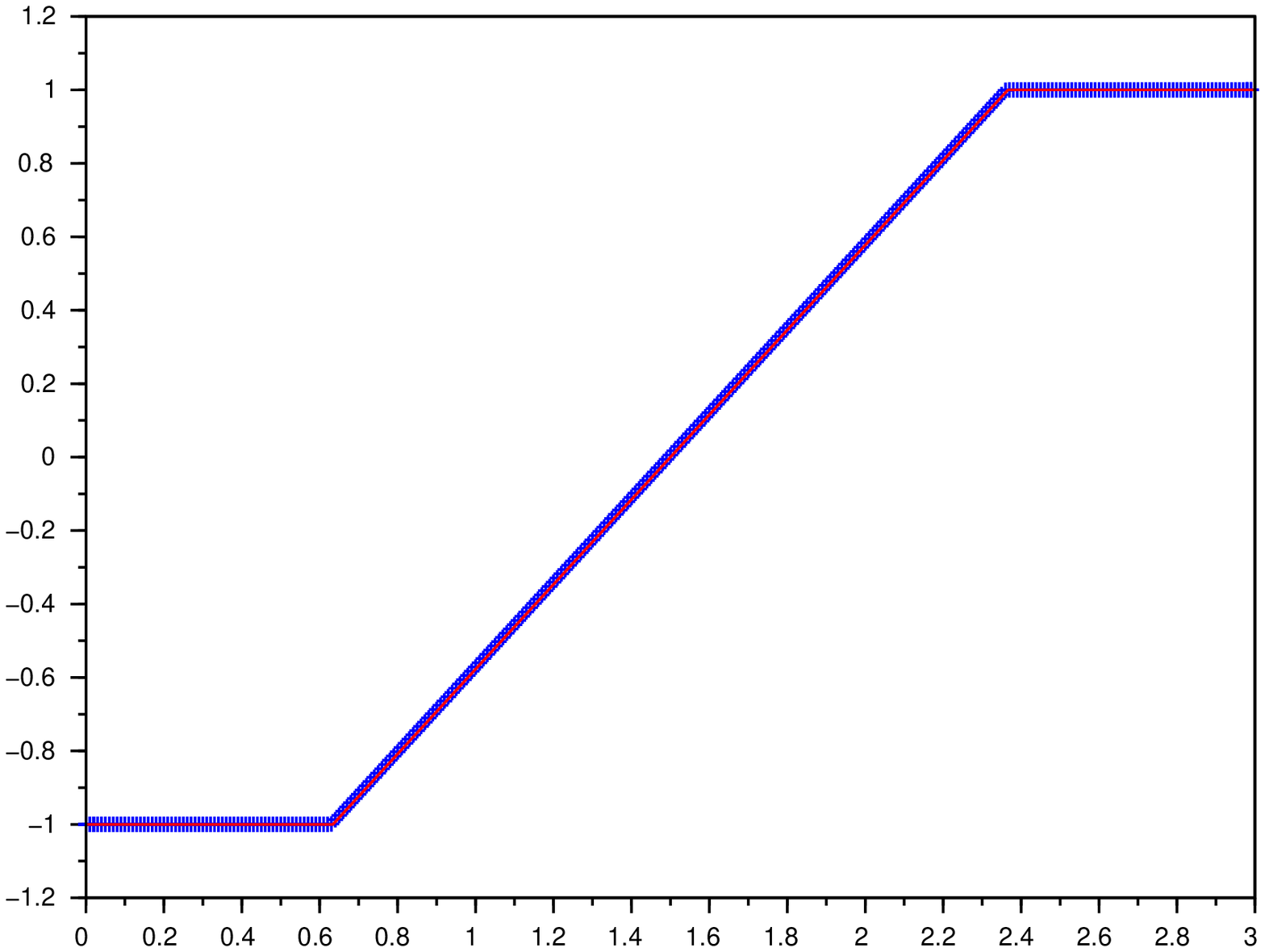} 
\end{center}
\item With $M=2$, $t_f=4$ (in the case $6M \leq t_f^2$) and for $T=1$, $T=0.5$, $T=0.1$ and $T=0.01$, we obtain:
\begin{center}
\includegraphics[scale=0.18]{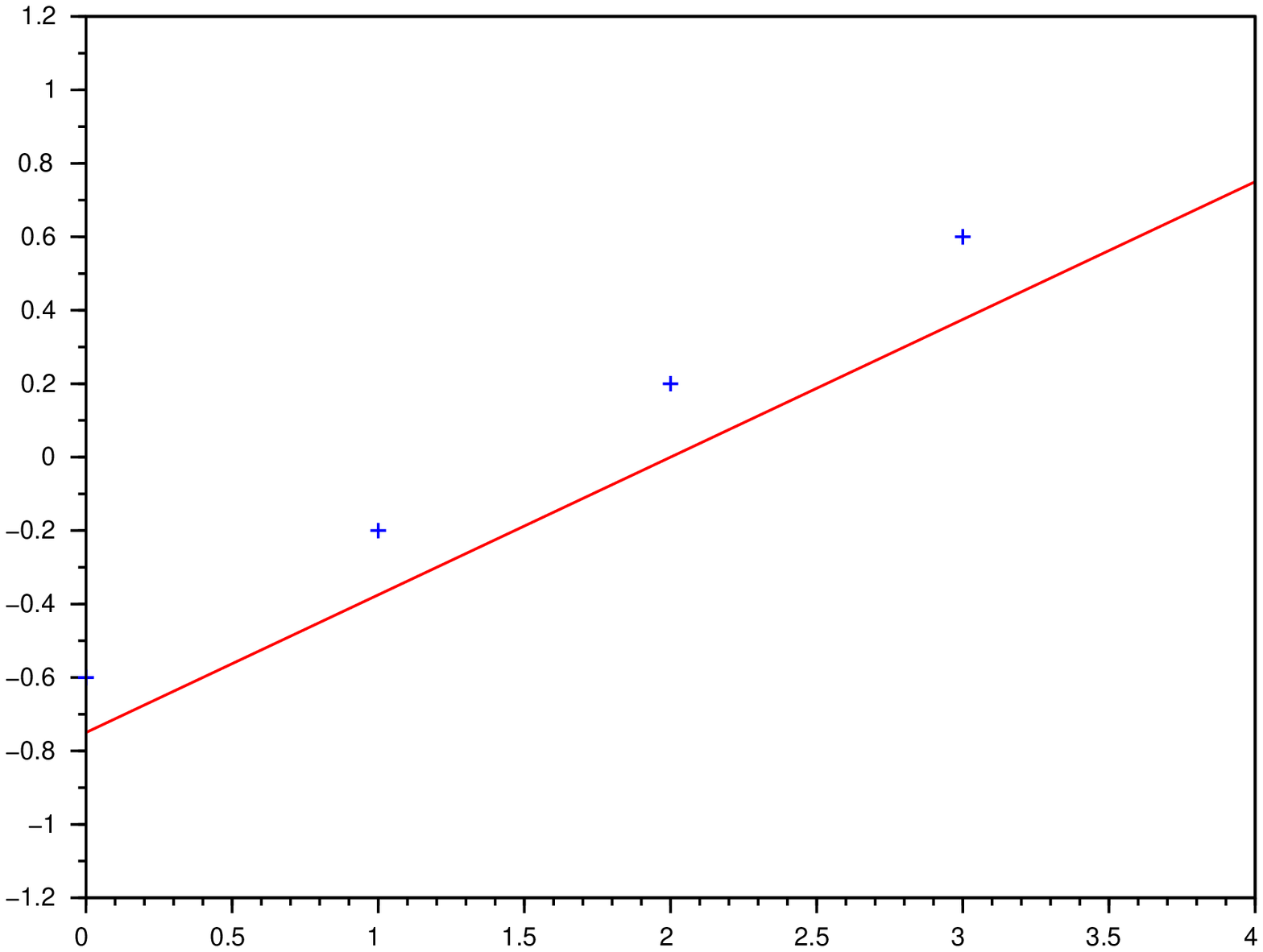}  \includegraphics[scale=0.18]{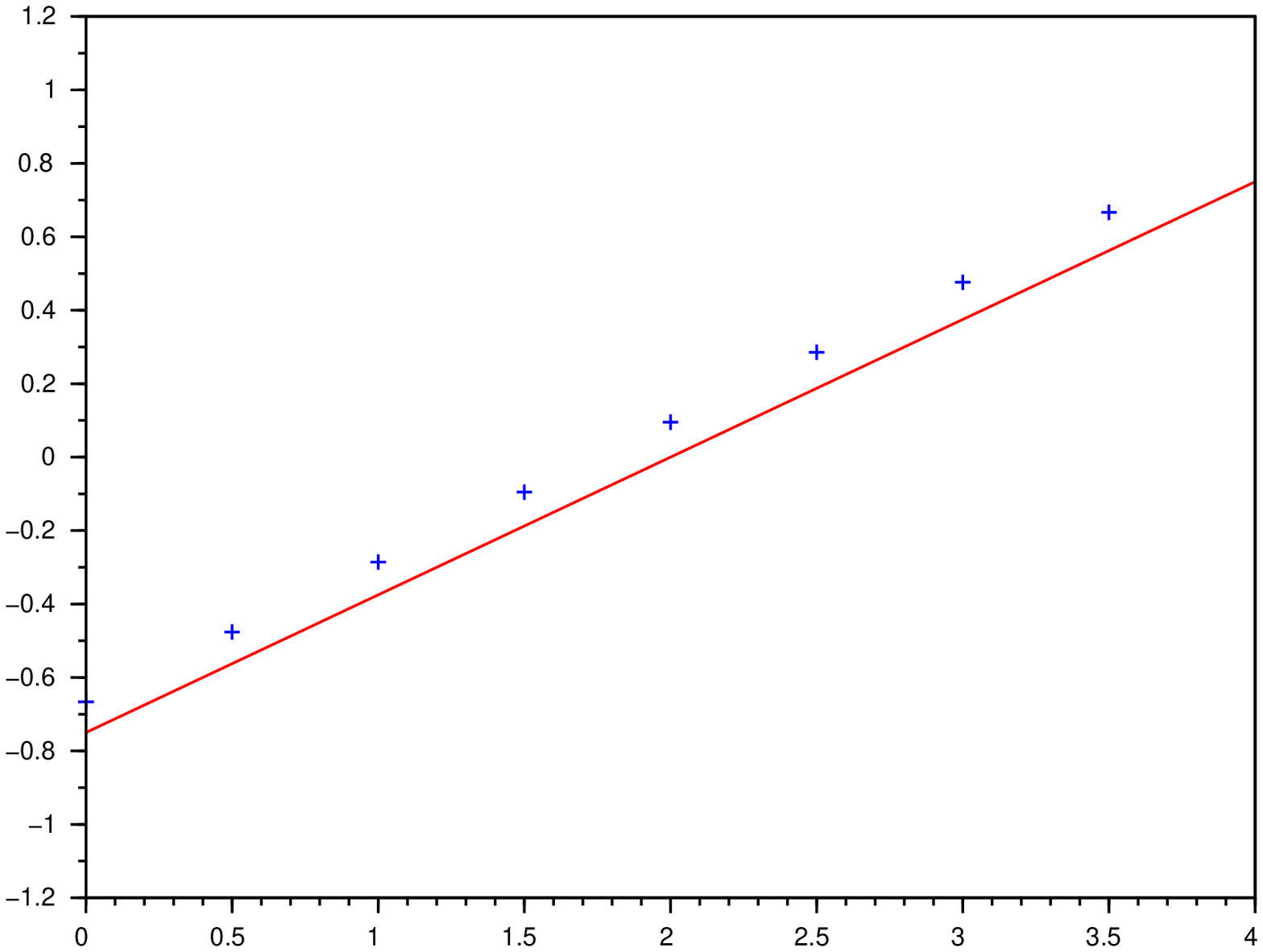} 
\includegraphics[scale=0.18]{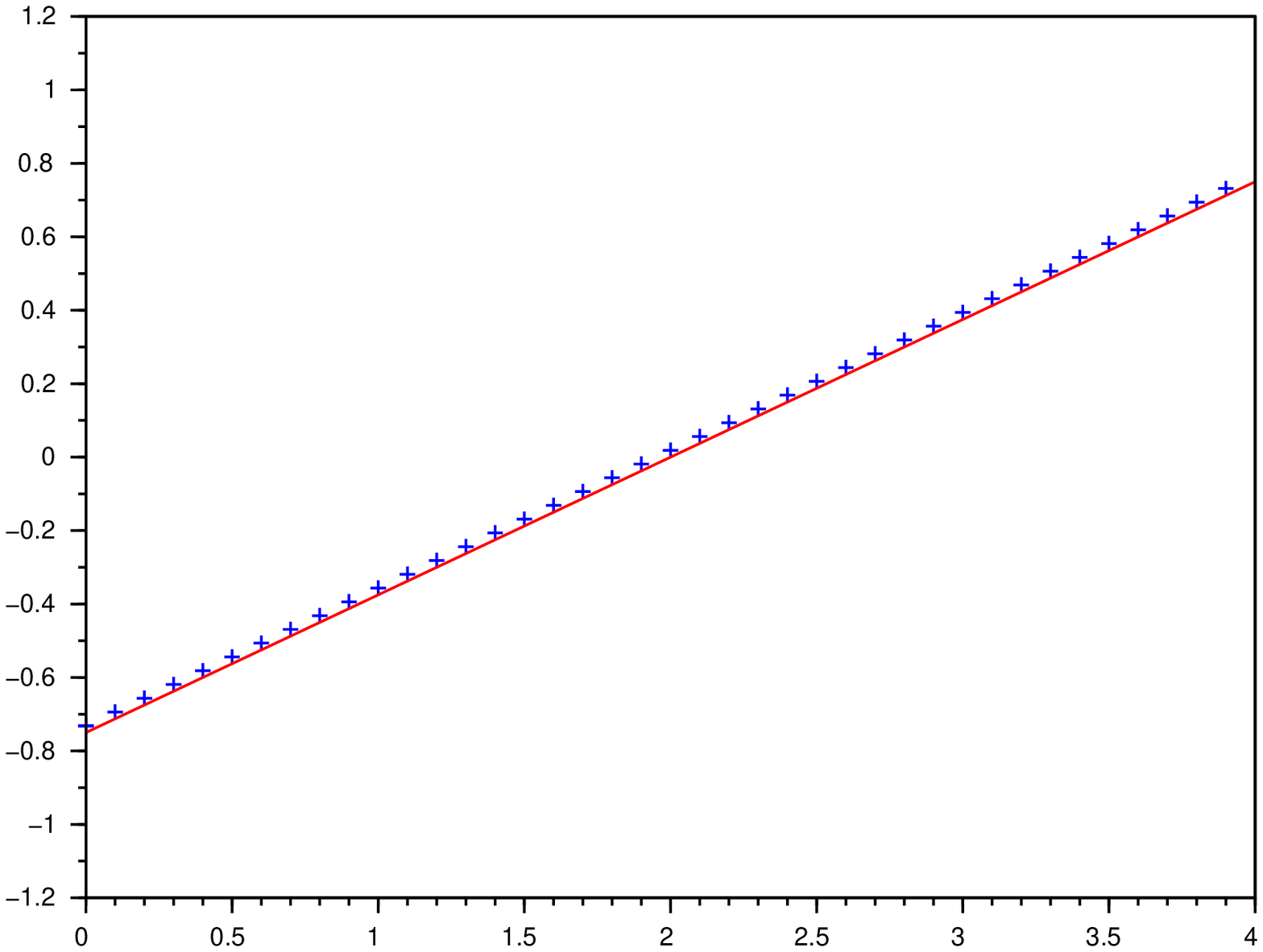}  \includegraphics[scale=0.18]{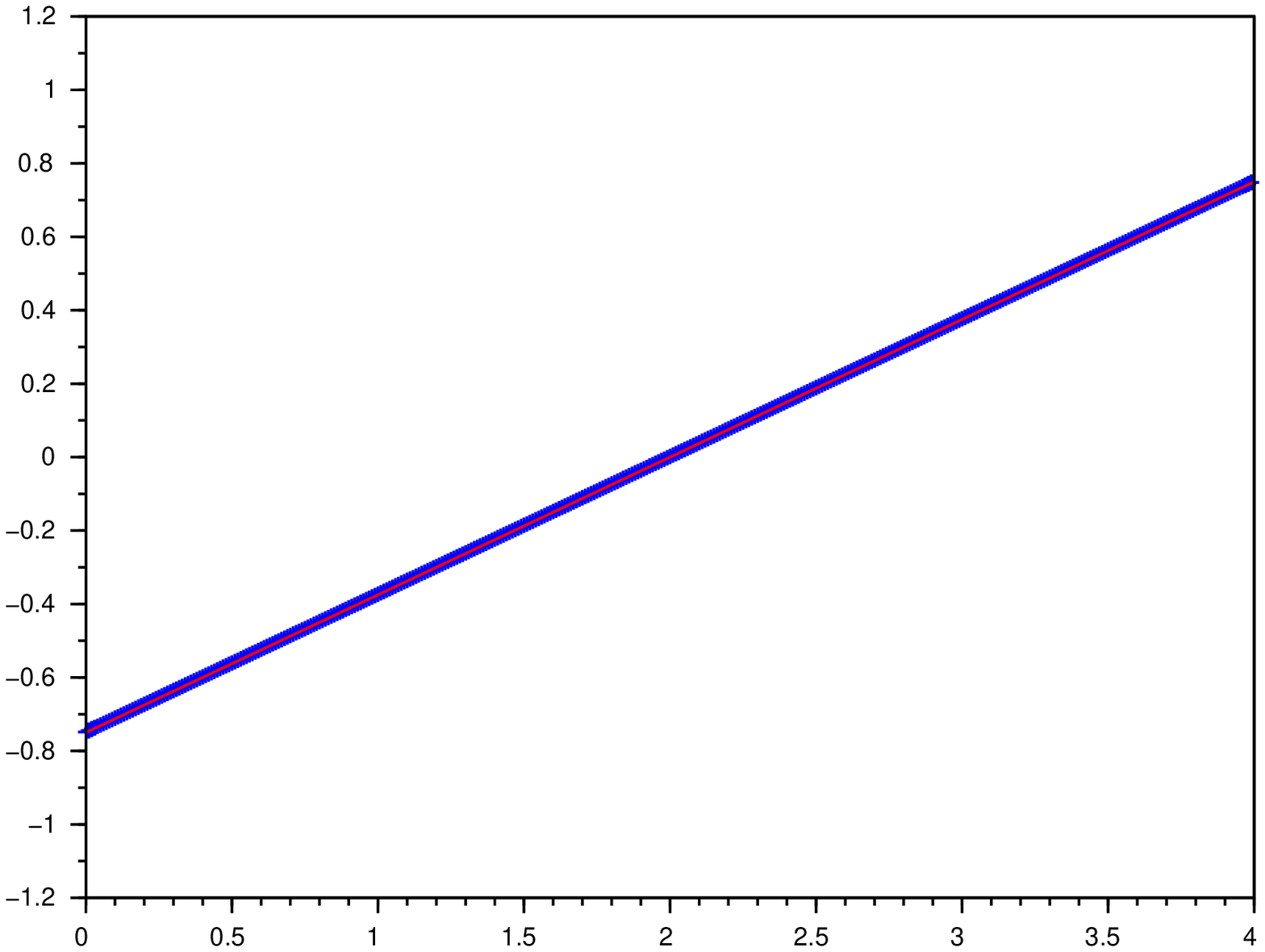} 
\end{center}

\end{itemize}

\begin{remark}
The previous numerical results naturally lead us to ask about the convergence of the optimal sampled-data control to the optimal (permanent) control when the sampling period $T$ tends to $0$. Actually, this natural question also emerges from the maximization condition~\eqref{secondconditionintro} that can be seen as an average of the weak maximization condition of the classical PMP, see Section~\ref{section12}. Indeed, note that the interval of average is smaller and smaller as the sampling period $T$ is reduced. Similarly, an important scientific perspective concerns the convergence of the optimal trajectory associated to a sampled-data control to the optimal trajectory associated to a permanent control. These important issues both constitute a forthcoming research project of the two authors of this note.
\end{remark}

\begin{remark}
Note that the above graphics only represent (by blue crosses) the discrete values $u(kT)$ of the sampled-data control $u$ at each controlling time $t=kT$. Let us provide some graphics representing the \textit{sample-and-hold} procedure consisting of ``freezing" the control at each controlling time $kT$ on the corresponding sampling time interval $[kT,(k+1)T)$. We fix $T = 0.5$ and we consider first $(M,t_f) = (2,3)$, then $(M,t_f) = (2,4)$. We obtain:
\begin{center}
\includegraphics[scale=0.18]{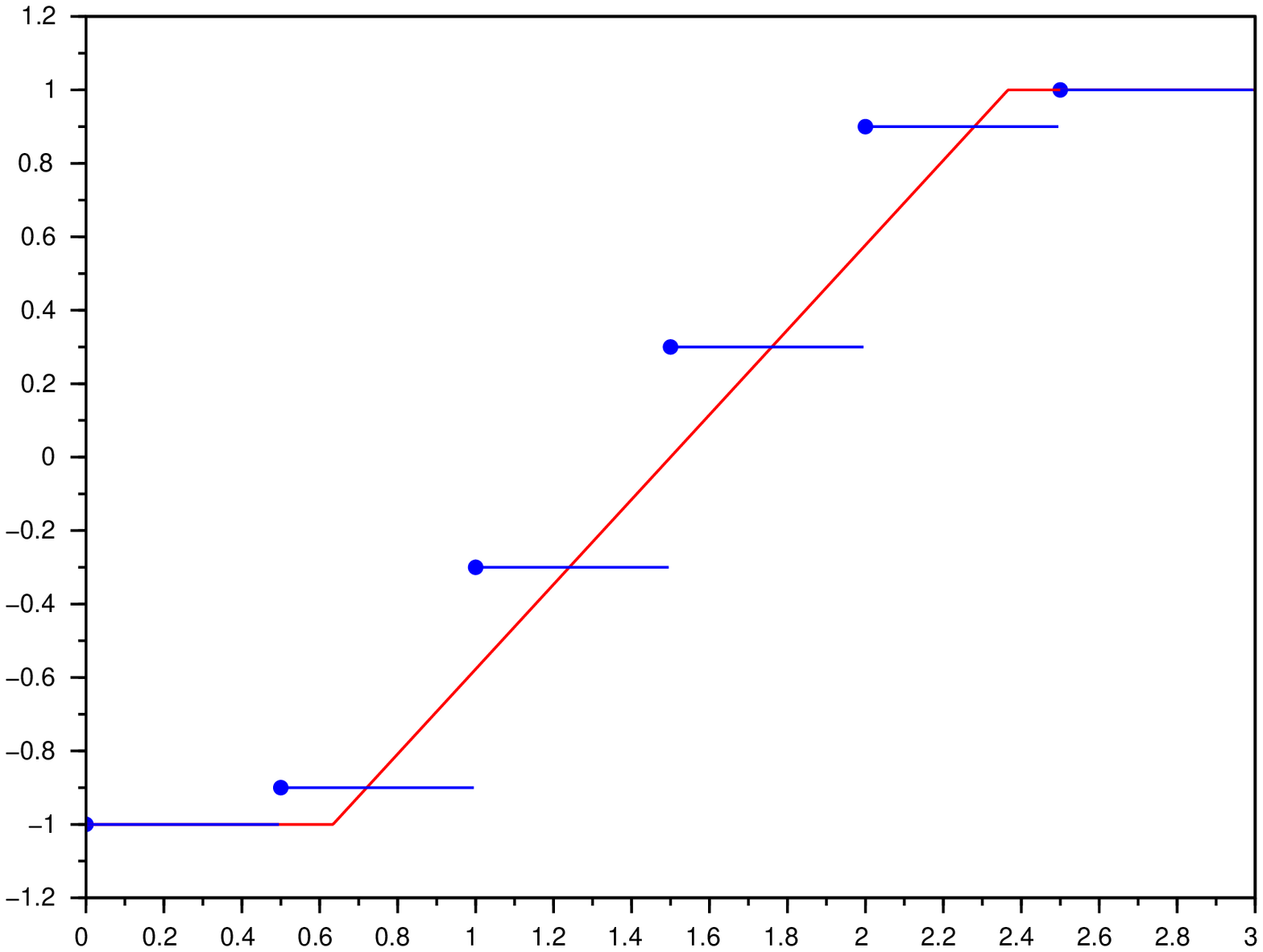}  \includegraphics[scale=0.18]{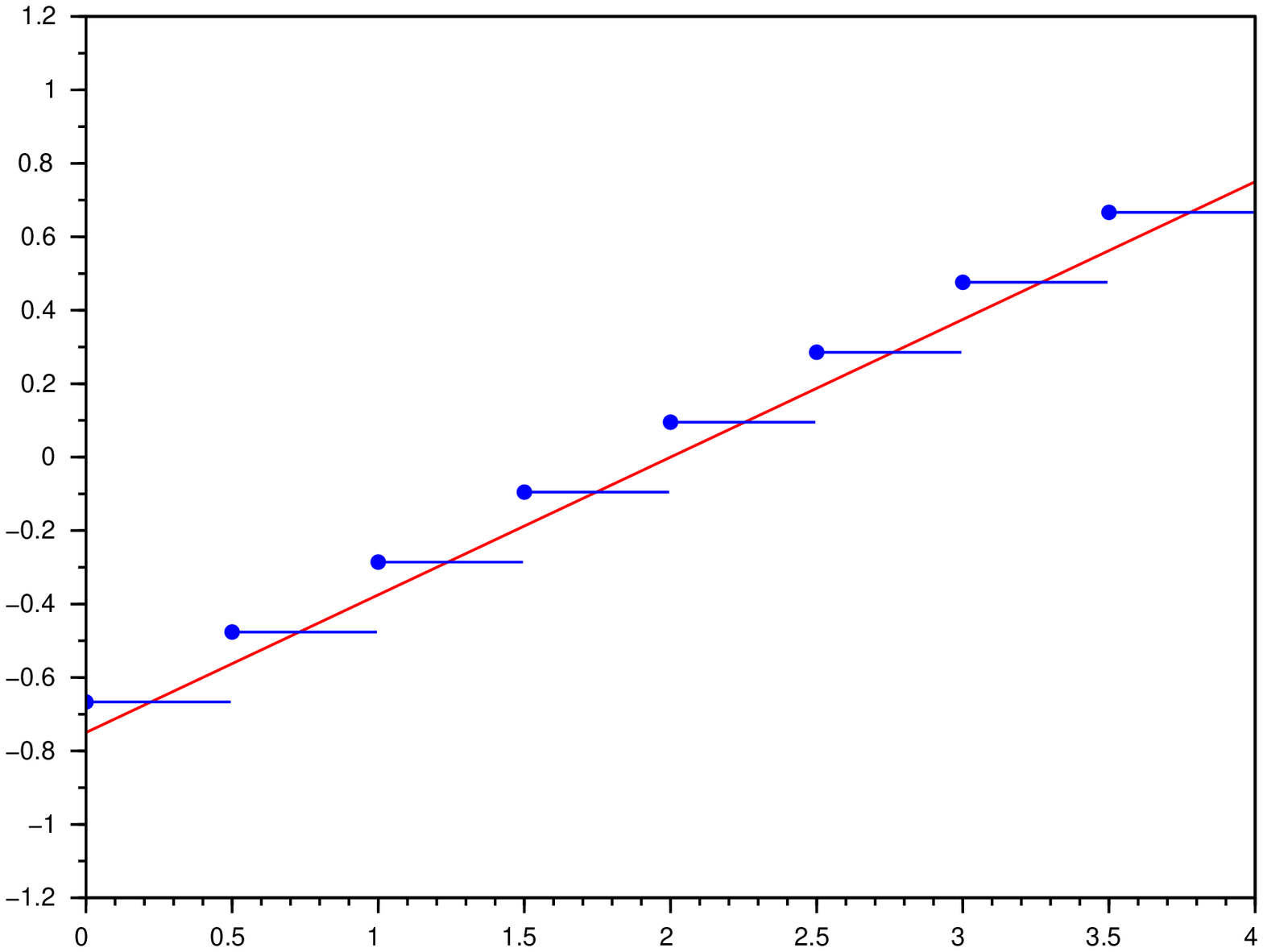} 
\end{center}
\end{remark}

\begin{ack}
The second author was partially supported by the Grant FA9550-14-1-0214 of the EOARD-AFOSR.
\end{ack}

%\bibliography{ifacconf}             % bib file to produce the bibliography
                                                     % with bibtex (preferred)

\end{document}